\documentstyle[12pt]{article}

\begin{document}
\begin{titlepage}
\begin{flushright}
math.QA/9907134
\end{flushright}
\vskip.3in

\begin{center}
{\Large \bf Level-one Highest Weight Representations of 
 $U_q[\widehat{gl(1|1)}]$ and Associated Vertex Operators  }
\vskip.3in
{\large Wen-Li Yang $^{a,b}$~~and Yao-Zhong Zhang $^{b}$}
\vskip.2in
{\em $~^a$ Institute of Modern Physics, Northwest University, 
Xian 710069 ,China \\
$~^b$ Department of Mathematics, University of Queensland, Brisbane,
     Qld 4072, Australia}
\end{center}

\vskip 2cm
\begin{center}
{\bf Abstract}
\end{center}

We study the level-one irreducible highest weight representations 
of  $U_q[\widehat{gl(1|1)}]$ and associated $q$-vertex operators.
We obtain the exchange relations satisfied by these vertex operators.
The characters and supercharacters associated with these irreducible
representations are calculated.


\end{titlepage}

%  Greek letters

\def\a{\alpha}
\def\b{\beta}
\def\d{\delta}
\def\e{\epsilon}
\def\ve{\varepsilon}
\def\g{\gamma}
\def\k{\kappa}
\def\l{\lambda}
\def\o{\omega}
\def\t{\theta}
\def\s{\sigma}
\def\D{\Delta}
\def\L{\Lambda}

\def\R{\overline{R}}
\def\G{{gl(N|N)}}
\def\hS{{\widehat{sl(N|N)}}}
\def\hG{{\widehat{gl(N|N)}}}
\def\R{{\cal R}}
\def\hR{{\hat{\cal R}}}
\def\C{{\bf C}}
\def\P{{\bf P}}
\def\Z2{{{\bf Z}_2}}
\def\Z{{\bf Z}}
\def\T{{\cal T}}
\def\H{{\cal H}}
\def\F{{\cal F}}
\def\V{\overline{V}}
\def\trho{{\tilde{\rho}}}
\def\tphi{{\tilde{\phi}}}
\def\tT{{\tilde{\cal T}}}
\def\uqsnh{{U_q[\widehat{sl(N|N)}]}}
\def\uqgnh{{U_q[\widehat{gl(N|N)}]}}
\def\uq1h{{U_q[\widehat{gl(1|1)}]}}
\def\uqg2h{{U_q[\widehat{gl(2|2)}]}}

% Shorthands for \begin{equation} and the like

\def\beq{\begin{equation}}
\def\eeq{\end{equation}}
\def\bea{\begin{eqnarray}}
\def\eea{\end{eqnarray}}
\def\ba{\begin{array}}
\def\ea{\end{array}}
\def\no{\nonumber}
\def\lt{\left}
\def\rt{\right}
\newcommand{\bq}{\begin{quote}}
\newcommand{\eq}{\end{quote}}

\newtheorem{Theorem}{Theorem}
\newtheorem{Definition}{Definition}
\newtheorem{Proposition}{Proposition}
\newtheorem{Lemma}{Lemma}
\newtheorem{Corollary}{Corollary}
\newcommand{\proof}[1]{{\bf Proof. }
        #1\begin{flushright}$\Box$\end{flushright}}

\newcommand{\sect}[1]{\setcounter{equation}{0}\section{#1}}
\renewcommand{\theequation}{\thesection.\arabic{equation}}

\sect{Introduction\label{intro}}

This paper is concerned with the level-one irreducible highest weight
representations and associated $q$-vertex operators of the simplest
quantum affine superalgebra $U_q[\widehat{gl(1|1)}]$.

Free bosonic realization of level-one representations and the 
corresponding $q$-vertex operators \cite{Fre92}
of quantum affine (bosonic) algebras have been investigated by a number
of groups (see e.g. \cite{Fre88,Ber89,Awa94,Jin95}). Such
kind of bosonization construction has been recently extended to the case of
type I quantum affine superalgebras $U_q[\widehat{sl(M|N)}],~M\neq N$
\cite{Kim97} and $U_q[\widehat{gl(N|N)}]$ \cite{Zha98a}. However,
the level-one irreducible highest weight representations and associated
$q$-vertex operators have been studied for $U_q[\widehat{sl(2|1)}]$
only \cite{Kim97,Yan99}. As is expected, the representation theory
of the super cases is much more complicated than that of the non-super cases.

It is well known by now that infinite dimensional irreducible highest weight 
representations and associated $q$-vertex operators play a very
powerful role in the algebraic analysis of massive integrable systems (see e.g.
\cite{Dav93,Jim94,Koy94,Yan99,Hou99}). Under some reasonable assumptions
on the physical space of states, this algebraic analysis method
\cite{Dav93,Jim94}
based on the infinite dimensional non-abelian quantum affine
(super)algebra symmetries enables one to compute the correlation
functions and form factors of massive (super) integrable systems in the form of
integral representations.

In this paper we study in details the level-one
irreducible highest weight representations of  $U_q[\widehat{gl(1|1)}]$
and associated vertex operators by using the free bosonic realization
given in \cite{Zha98a}.  
We calculate  the exchange relations satisfied by the vertex operators,
and compute the characters and supercharacters associated with
these irreducible representations .

\sect{Bosonization of $\uq1h$ at Level-One}

\subsection{Drinfeld basis of  $\uq1h$}

The simple
roots for $\widehat{gl(1|1)}$  are
$\a_0=\d-\ve_1+\ve_{2}$, $\a_1=\ve_1-\ve_{2}$ 
with $\d,~\{\ve_1,~\ve_2\}$ satisfying
\beq
(\d,\d)=(\d,\ve_k)=0,~~~~(\ve_k,\ve_{k'})=(-1)^{k+1}\d_{kk'}, ~~
k,k'=1,2.
\eeq 
The generalized symmetric Cartan matrix of 
$\widehat{gl(1|1)}$  is degenerate. For the reason which will become clear
later in the construction of the vertex operator of $\uq1h$, we
extend the Cartan subalgebra \cite{Zha98a}
by adding to it the element $\a_2=\ve_1+\ve_2$.
The enlarged Cartan matrix of $\widehat{gl(1|1)}$ has elements
$a_{ij}=(\a_i,\a_j)$, $i,j=0,1,2$, so that the Cartan matrix
$(a_{ij}),~i,j=1,2$ of $gl(1|1)$ is invertible.
Denote by $\H$ the extended Cartan
subalgebra and by $\H^*$ the dual of $\H$.
Let $\{h_0,h_1,h_{2},d\}$ be a basis of $\H$, where  $d$
is the usual derivation operator.
Let $\{\L_0,\L_1,\L_{2},\delta\}$
be the dual basis with $\L_j$ being fundamental weights.
Explicitly \cite{Zha98a}
\bea
&&\L_{2}=\frac{\ve_1-\ve_2}{2}, ~~~~
\L_1=\L_0+\frac{\ve_1+\ve_2}{2},~~~~\L_0.
\eea
The quantum affine superalgebra $\uq1h$ is a quantum (or $q$-)
deformation of the universal enveloping algebra of $\widehat{gl(1|1)}$
 and is
generated by the Chevalley generators $\{e_i,\;f_i\;q^{h_j},\;d|
i=0,1,~j=0,1,2\}$.
The $\Z_2$-grading of the Chevalley generators is $[e_i]=[f_i]=1,~
i=0,1$ and zero
otherwise. The defining relations are
\bea
&&hh'=h'h,~~~~~~\forall h\in \H,\no\\
&&q^{h_j}e_iq^{-h_j}=q^{a_{ij}}e_i,~~~~[d, e_i]=\d_{i0}e_i,\no\\
&&q^{h_j}f_iq^{-h_j}=q^{-a_{ij}}f_i,~~~~[d,f_i]=-\d_{i0}f_i,\no\\
&&[e_i,f_{i'}]=\d_{ii'}\frac{q^{h_i}-q^{-h_i}}{q-q^{-1}},\no\\
&&[e_i,e_{i'}]=[f_i,f_{i'}]=0,~~~~{\rm for}~~a_{ii'}=0,\no\\
&&[[e_0,e_1]_{q^{-1}}, [e_0,e_{1}]_q]=0,~~
[[f_0,f_1]_{q^{-1}}, [f_0,f_{1}]_q]=0.
\eea
Here and throughout, $[a,b]_x\equiv ab-(-1)^{[a][b]}x ba$ 
and $[a,b]\equiv [a,b]_1$.

$\uq1h$ is a $\Z_2$-graded quasi-triangular Hopf algebra endowed with
the following coproduct $\D$, counit $\e$ and antipode $S$:
\bea
\D(h)&=&h\otimes 1+1\otimes h,\no\\
\D(e_i)&=&e_i\otimes 1+q^{h_i}\otimes e_i,~~~~
\D(f_i)=f_i\otimes q^{-h_i}+1\otimes f_i,\no\\
\e(e_i)&=&\e(f_i)=\e(h)=0,\no\\
S(e_i)&=&-q^{-h_i}e_i,~~~~S(f_i)=-f_iq^{h_i},~~~~S(h)=-h,\label{e-s} 
\eea
where $i=0,1$ and $h\in \H$.

$\uq1h$ can also be realized in terms of the Drinfeld generators
\cite{Dri88} $\{X^{\pm}_m,\; H^j_n,\; q^{\pm H^j_0}$, $c,\; d |
m\in {\bf Z},\; n\in{\bf Z}-\{0\},\; 
j=1,2\}$.  The $\Z_2$-grading of the Drinfeld generators is
given by $[X^{\pm}_m]=1$, for $m\in{\bf Z}$ and 
$[H^j_n]=[H^j_0]=[c]=[d]=0$ for all $j=1,2,\; n\in{\bf Z}-\{0\}$.
The relations 
satisfied by the Drinfeld generators read \cite{Yam96,Zha97,Zha98a}
\bea
&&[c,a]=[d,H^j_0]=[H^j_0, H^{j'}_n]=0,~~~~\forall a\in\uq1h\no\\
&&q^{H^j_0}X^{\pm}_nq^{-H^j_0}=q^{\pm a_{1j}}X^{\pm}_n,\no\\
&&[d,X^{\pm}_n]=nX^{\pm}_n,~~~[d,H^j_n]=nH^j_n,\no\\
&&[H^j_n, H^{j'}_m]=\d_{n+m,0}\frac{[a_{jj'}n]_q[nc]_q}{n},\no\\
&&[H^j_n,
   X^{\pm}_m]=\pm\frac{[a_{1j}n]_q}{n}X^{\pm}_{n+m}q^{\mp|n|c/2},\no\\
&&[X^{+}_n, X^{-}_m]=\frac{1}{q-q^{-1}}\lt(q^{\frac{c}{2}(n-m)}
  \psi^{+,1}_{n+m}-q^{-\frac{c}{2}(n-m)}\psi^{-,1}_{n+m}\rt),\no\\
&&[X^{\pm}_n, X^{\pm}_m]=0~~.
%&&[X^{\pm,i}_{n+1}, X^{\pm,i'}_m]_{q^{\pm a_{ii'}}}
%  -[X^{\pm,i'}_{m+1}, X^{\pm,i}_n]_{q^{\pm a_{ii'}}}=0,\no\\
%&&[[X^{\pm,l}_m,X^{\pm,l-1}_{m'}]_{q^{(-1)^l}},[X^{\pm,l}_n,
%  X^{\pm,l+1}_{n'}]_{q^{(-1)^{l+1}}}]\no\\
%&&~~~+[[X^{\pm,l}_n,X^{\pm,l-1}_{m'}]_{q^{(-1)^l}},[X^{\pm,l}_m,
%  X^{\pm,l+1}_{n'}]_{q^{(-1)^{l+1}}}]=0, ~~~l=2,\cdots,2N-2.
\label{drinfeld}
\eea
where $[x]_q=(q^x-q^{-x})/(q-q^{-1})$ and
 $\psi^{\pm,j}_{n}$ are related to $H^j_{\pm n}$ by relations
\beq
\sum_{n\in{\bf Z}}\psi^{\pm,j}_{n}z^{- n}=q^{\pm H^j_0}\exp\lt(
  \pm(q-q^{-1})\sum_{n>0}H^j_{\pm n}z^{\mp n}\rt).
\eeq
The Chevalley generators are related to the Drinfeld generators by the
formulae
\bea
&&h_i=H_0^i,~~~e_1=X^{+}_0,~~~f_1=X^{-}_0,~~~
e_0=X^-_{1}q^{-H^1_0},~~~f_0=-q^{H^1_0}X^+_{-1},\no\\
&&h_{2N}=H^{2N}_0,~~~h_0=c-H^1_0,
\eea
\subsection{Level-one free bosonic realization}
In this subsection, we briefly review the bosonization of 
 $\uq1h$ at level one \cite{Zha98a}.

Let us introduce bosonic oscillators $\{a^j_n,\;c_n,\;Q_{a^j},\;
Q_{c}|n\in{\bf Z}, j=1,2,\}$ which
satisfy the commutation relations
\bea
&&[a^i_n, a^{j}_m]=(-1)^{i+1}\d_{ij}\d_{m+n,0}\frac{[n]^2_q}{n},~~~~~
  [a^i_0, Q_{a^{j}}]=\d_{ij},~~~~i,j=1,2,\no\\
&&[c_n, c_m]=\d_{n+m,0}\frac{[n]_q^2}{n},~~~~~
  [c_0, Q_{c}]=1.\label{oscilators}
\eea
The remaining commutation relations are zero. Corresponding to these
bosonic oscillators are the $q$-deformed free bosonic currents
\bea
&&H^j(z;\k)=Q_{A^j}+A^j_0\ln z   
   -\sum_{n\neq 0}\frac{A^j_n}{[n]_q}q^{\k |n|}z^{-n},\no\\
&&c(z)=Q_{c}+c_0\ln z-\sum_{n\neq 0}\frac{c_n}{[n]_q}  
   z^{-n},\no\\
&&H^j_\pm(z)=
   \pm(q-q^{-1})\sum_{n>0}A^j_{\pm n}z^{\mp n}\pm A^i_0\ln q,
\eea
where 
\bea
&&A^1_n=a^1_n+a^{2}_n,~~~~
  A^{2}_n=\frac{q^n+q^{-n}}{2}(a^1_n-a^2_n),\no\\
&&Q_{A^1}=Q_{a^1}-Q_{a^{2}},~~~~
  Q_{A^{2}}=Q_{a^1}+Q_{a^2}.\label{A-a}
\eea
\noindent We introduce the Drinfeld currents or generating 
functions
\bea
X^{\pm}(z)=\sum_{n\in\Z}X^{\pm}_nz^{-n-1}~,~~
\psi^{\pm,j}(z)=\sum_{n\in\Z}\psi^{\pm,j}_nz^{-n}~, ~~j=1,2,
\eea
\noindent and the q-differential operator defined by 
$\partial_zf(z)=\frac{f(qz)-f(q^{-1}z)}{(q-q^{-1})z}$.
Then ,   
\begin{Theorem}\label{free-boson}(\cite{Zha98a}):
The Drinfeld generators of $\uq1h$  at level one are
realized by the free boson fields as
\bea
&&c=1,\\
&&\psi^{\pm,j}(z)=e^{H^j_\pm(z)},~~~~j=1,2,\\
&&X^{+}(z)=:e^{ H^1(z;-\frac{1}{2})}\;e^{c(z)}:,\\
&&X^{-}(z)=:e^{ -H^1(z;\frac{1}{2})}\;\partial_z\{e^{-c(z)}\}:.
\eea
\end{Theorem}

\subsection{Bosonization of level-one vertex operators}

We consider the evaluation representation $V_z$ of $\uq1h$,
where $V$ is a two-dimensional graded
vector space with basis vectors $\{v_1,v_2\}$. The
$\Z_2$-grading of the basis vectors is chosen to be $[v_j]=\frac{(-1)^j
+1}{2}$. Let $e_{j,j'}$ be the $2\times 2$ matrices satisfying
$e_{i,j}v_k=\d_{jk}v_i$. Let $V^{*S}_z$ be the dual module of $V_z$
defined by $\pi_{V^{*S}}(a)=\pi_V(S(a))^{st},~ \forall a\in 
U_q[\widehat{gl(1|1)}]$, where $st$ is the supertransposition operation.

In the homogeneous gradation,
the Drinfeld generators are represented
on $V_z$ by \cite{Zha98a}
\bea
H^1_m&=&\frac{[m]_q}{m}z^m
  (e_{1,1}+e_{2,2}),~H^{2}_m=-z^m\frac{[2m]_q}{m}q^me_{2,2},
~H^{2}_0=-2e_{2,2}\no\\
H^1_0&=&e_{1,1}+e_{2,2},~X^{+}_m=(qz)^me_{1,2},~
  X^{-}_m=(qz)^me_{2,1},
\eea
and on $V^{*S}_z$ by
\bea
H^1_m&=&-\frac{[m]_q}{m}z^m (e_{1,1}+e_{2,2}),~
H^{2}_m=z^m\frac{[2m]_q}{m}q^{-m}e_{2,2}~,
~~H^{2}_0=2e_{2,2},\no\\
H^1_0&=&-e_{1,1}-e_{2,2},~~
X^{+}_m=q^{-1}(q^{-1}z)^me_{2,1},~
  X^{-}_m=-q(q^{-1}z)^me_{1,2}.
\eea

Now,  let 
$V(\l)$ be the highest weight $\uq1h$-module with the highest weight
$\l$. Consider the following intertwiners of $\uq1h$-modules
\cite{Jim94}:
\bea
\Phi^{\mu V}_\l(z)&:&~~ V(\l)\longrightarrow V(\mu)\otimes V_z,
     \label{Phi}\\
\Phi^{\mu V^*}_\l(z)&:&~~ V(\l)\longrightarrow V(\mu)\otimes V^{*S}_z,
     \label{Phi*}\\
\Psi^{V\mu}_\l(z)&:&~~ V(\l)\longrightarrow V_z\otimes V(\mu),
     \label{Psi}\\
\Psi^{V^*\mu}_\l(z)&:&~~ V(\l)\longrightarrow V^{*S}_z\otimes V(\mu).
     \label{Psi*}
\eea
They are intertwiners in the sense that for any $x\in \uq1h$
\beq
\Xi(z)\cdot x=\D(x)\cdot\Xi(z),~~~~\Xi(z)=\Phi^{\mu V}_\l(z),~
   \Phi^{\mu V^*}_\l(z),~\Psi^{V\mu}_\l(z),~\Psi^{V^*\mu}_\l(z).
   \label{intertwiner1}
\eeq
These intertwiners are even operators, that is their gradings are
$[\Phi^{\mu V}_\l(z)] = [\Phi^{\mu V^*}_\l(z)] = [\Psi^{V\mu}_\l(z)]
= [\Psi^{V^*\mu}_\l(z)] = 0$. According to \cite{Jim94},
$\Phi^{\mu V}_\l(z)~\lt(\Phi^{\mu V^*}_\l(z)\rt)$ is called type I
(dual) vertex operator and
$\Psi^{V\mu}_\l(z)~\lt(\Psi^{V^*\mu}_\l(z)\rt)$ type II (dual) vertex
operator. 
The vertex operators can be expanded in terms of the basis \cite{Jim94}
\bea
&&\Phi^{\mu V}_\l(z)=\sum_{j=1}^{2}\,\Phi^{\mu V}_{\l,j}(z)
   \otimes v_j, ~~~~
\Phi^{\mu V^*}_\l(z)=\sum_{j=1}^{2}\,\Phi^{\mu V^*}_{\l,j}(z)
   \otimes v^*_j,\no\\ 
&&\Psi^{V\mu}_\l(z)=\sum_{j=1}^{2}\,v_j\otimes \Psi^{V\mu}_{\l,j}(z),~~~~
\Psi^{V^*\mu}_\l(z)=\sum_{j=1}^{2}\,v^*_j\otimes \Psi^{V^*\mu}_{\l,j}(z).
\eea

The intertwining operators which satisfy (\ref{intertwiner1}) for 
any $x\in U_q[\widehat{sl(1|1)}]$ have been constructed in \cite{Zha98a}.
We extend the construction to $\uq1h$ by requiring that the vertex 
operators also obey (\ref{intertwiner1}) for the element $x=H^{2}_m$,
which extends  $U_q[\widehat{sl(1|1)}]$ to $\uq1h$.

Define 
the even operators
\bea
&&\phi(z)=\sum_{j=1}^{2}\,\phi_j(z)
   \otimes v_j, ~~~~
\phi^*(z)=\sum_{j=1}^{2}\,\phi^*_j(z)
   \otimes v^*_j,\no\\
&&\psi(z)=\sum_{j=1}^{2}\,v_j\otimes \psi_j(z),~~~~
\psi^*(z)=\sum_{j=1}^{2}\,v^*_j\otimes \psi^*_j(z).
\eea
Assuming that $\phi(z)$, $\phi^*(z)$, $\psi(z)$ and $\psi^*(z)$ 
satisfy (\ref{intertwiner1}) for any $x\in\uq1h$ and  by using
the results in \cite{Zha98a} for $U_q[\widehat{sl(1|1)}]$, we find
\bea
&&\phi_{2}(z)=:e^{-H^{*,1}(qz;\frac{1}{2})+H^1(zq^2;\frac{1}{2})}
e^{c(qz)}:e^{-\sqrt{-1}\pi c_0},~~~~
-\phi_1(z)=[\phi_{2}(z)~,~f_1]_{q^{-1}},\no\\
&&\phi^*_1(z)=:e^{H^{*,1}(qz;\frac{1}{2})}:e^{\sqrt{-1}\pi c_0},~~~~
q\phi^*_{2}(z)=[\phi^*_1(z)~,f_1]_{q},\no\\
&&\psi_1(z)=:e^{-H^{*,1}(qz;-\frac{1}{2})}:
  e^{-\sqrt{-1}\pi c_0}~~,~~
  \psi_{2}(z)=[\psi_1(z)~,e_1]_{q},\no\\
&&\psi^*_{2}(z)=:e^{H^{*,1}(qz;-\frac{1}{2})-H^1(z;-\frac{1}{2})}
\partial_z\{e^{-c(qz)}\}:e^{\sqrt{-1}\pi c_0},\no\\~            
&&q^{-1}\psi^*_1(z)=[\psi^*_{2}(z)~,~e_1]_{q^{-1}}.
\label{Vertex-operator}
\eea
\noindent where 
\bea
&&H^{*,1}(z;\k)=Q^*_{A^1}+A^{*1}_0\ln z-\sum_{n\neq
0}\frac{A^{*1}_n}{[n]_q} q^{k|n|}z^{-n},\\
&&A^{*1}_n=\frac{1}{q^n+q^{-n}}A^{2}_n,~~
A^{*2}_n=A^{1}_n,~~n\neq 0,\\
&&A^{*1}_0=\frac{1}{2}A^2_0,~~A^{*2}_0=\frac{1}{2}A^1_0,~~
  Q^*_{A^1}=\frac{1}{2}Q_{A^2},~~
  Q^*_{A^2}=\frac{1}{2}Q_{A^1}.  
\eea
\noindent Since $\phi(z)$, $\phi^*(z)$, $\psi(z)$ and $\psi^*(z)$ 
satisfy the same intertwining relations as $\Phi^{\mu V}_\l(z)$, $
\Phi^{\mu V^*}_\l(z)$, $\Psi^{V\mu}_\l(z)$ and $\Psi^{V^*\mu}_\l(z)$
respectively, we have 
\begin{Proposition}:
The vertex operators $\Phi^{\mu V}_\l(z),\;
\Phi^{\mu V^*}_\l(z),\; \Psi^{V\mu}_\l(z)$ and $\Psi^{V^*\mu}_\l(z)$, if
they exist, have  the same bosonizaion as the operators
$\phi(z),\;\phi^*(z),\; \psi(z)$ and $\psi^*(z)$, respectively. 
\end{Proposition}

\section{Exchange Relations of the Bosonized Vertex Operators}
In this section, we derive the exchange relations of the type I and
type II bosonized vertex operators of $\uq1h$. As expected, these 
vertex operators satisfy the graded Faddeev-Zamolodchikov algebra.

Let $R(z) \in End(V\otimes V)$ be the R-matrix of  $\uq1h$, 
defined by 
\bea
R(z)(v_i\otimes v_j)=\sum_{k,l=1}^{2}R^{ij}_{kl}(z)v_k\otimes v_l,~~~~
\forall v_i,v_j,v_k,v_l\in V,
\eea
\noindent where 
\begin{eqnarray*}
& &R^{1,1}_{1,1}(z)=1,~~R^{2,2}_{2,2}(z)=
\frac{zq^{-1}-q}{zq-q^{-1}},~~
R^{21}_{12}(z)=\frac{q-q^{-1}}{zq-q^{-1}}~,\\
& &R^{12}_{12}(z)=R^{21}_{21}(z)=\frac{z-1}{zq-q^{-1}},~~
R^{12}_{21}(z)=\frac{(q-q^{-1})z}{zq-q^{-1}}~,
~~R^{ij}_{kl}(z)=0,~~{\rm otherwise}.
\end{eqnarray*}
\noindent  
The R-matrix satisfies the graded Yang-Baxter equation  on 
 $V\otimes V\otimes V$
\begin{eqnarray*}
R_{12}(z)R_{13}(zw)R_{23}(w)=R_{23}(w)R_{13}(zw)R_{12}(z),
\end{eqnarray*}
\noindent moreover enjoys : (i) initial condition, $R(1)=P$ with 
$P$ being the graded permutation operator; (ii) unitarity condition, 
$R_{12}(\frac{z}{w})R_{21}(\frac{w}{z})=1$ , where $R_{21}(z)=P
R_{12}(z)P$; and (iii) crossing-unitarity,
\begin{eqnarray*}
R^{-1,st_1}(z)R(z)^{st_1}=\frac{(z-1)^2}{(q^{-1}z-q)(zq-q^{-1})}.
\end{eqnarray*}
\noindent The various supertranspositions of the R-matrix are given by 
\begin{eqnarray*}
& &(R^{st_1}(z))^{kl}_{ij}=R^{il}_{kj}(z)(-1)^{[i]([i]+[k])},~~~~
(R^{st_2}(z))^{kl}_{ij}=R^{kj}_{il}(z)(-1)^{[j]([l]+[j])},\\
& &(R^{st_{12}}(z))^{kl}_{ij}=R^{ij}_{kl}(z)
(-1)^{([i]+[j])([i]+[j]+[k]+[l])}=R^{ij}_{kl}(z).
\end{eqnarray*}

Now we calculate the exchange relations of the type I and type 
II bosonic vertex operators of $\uq1h$ in (\ref{Vertex-operator}).
Define 
\begin{eqnarray*}
\oint dzf(z)=Res(f)=f_{-1},~~~~{\rm for~ a~ formal~ series~ function} 
~ f(z)=\sum_{n\in \Z}f_nz^n.
\end{eqnarray*}
\noindent Then, the Chevalley generators of $\uq1h$ can be expressed by
the integrals
\begin{eqnarray*}
e_1=\oint dz X^{+}(z),~~~~f_1=\oint dz X^{-}(z).
\end{eqnarray*}
\noindent One can also get the integral expression of the bosonic vertex
operators $\phi(z)$, $\phi^*(z)$, $\psi(z)$ and $\phi^*(z)$ from  
(\ref{Vertex-operator}). Using these integral expressions , we arrive at  
\begin{Proposition}
: The bosonic vertex operators defined in (\ref{Vertex-operator})
satisfy the graded Faddeev-Zamolodchikov algebra
\bea
&&\phi_j(z_2)\phi_i(z_1)=
\sum_{k,l=1}^{2}R^{kl}_{ij}(\frac{z_1}{z_2})
\phi_k(z_1)\phi_l(z_2)(-1)^{[i][j]},\label{ZF1}\\
&&\psi^*_i(z_1)\psi^*_j(z_2)=
\sum_{k,l=1}^{2}R^{ij}_{kl}(\frac{z_1}{z_2})
\psi^*_l(z_2)\psi^*_k(z_1)(-1)^{[i][j]},\\
&&\psi^*_i(z_1)\phi_j(z_2)=
\phi_j(z_2)\psi^*_i(z_1)(-1)^{[i][j]},\label{ZF2}
\eea
\noindent and the following invertibility relations
\begin{eqnarray*}
\phi_i(z)\phi^*_j(z)=-q\d_{ij}~id.
\end{eqnarray*}

\end{Proposition}
In the derivation of this proposition the fact that $R^{kl}_{ij}(z)
(-1)^{[k][l]}=R^{kl}_{ij}(z)(-1)^{[i][j]}$ is helpful.

\sect{Irreducible Highest Weight $\uq1h$-modules at Level-One }

In this section we study in details the irreducible $\uq1h$-module structure 
in the Fock  space.

We begin by defining the Fock module. Denote by
$F_{\l_1,\l_2;\l_3}$ the bosonic Fock
spaces generated by $a_{-m}^i,c_{-m} (m>0)$
 over the vector $|\l_1,\l_2;\l_3>$:
\begin{eqnarray*}
F_{\l_1,\l_2;\l_3}=
{\bf C}[a^1_{-1},a^2_{-1}, a^1_{-2},a^2_{-2}, \cdots;c_{-1},c_{-2},\cdots]
|\l_1,\l_2;\l_3>,
\end{eqnarray*}
\noindent where 
\begin{eqnarray*}
|\l_1,\l_2;\l_3>=
e^{\sum_{i=1}^{2}\l_iQ_{a^i}+\l_3Q_{c}}|0>.
\end{eqnarray*}
\noindent The vacuum vector $|0>$ is defined by 
$a^i_m|0>=c_m|0>=0$ for $i=1,2$ and $m\geq 0$. Obviously, 
\begin{eqnarray*}
a^i_m|\l_1,\l_2;\l_3>=0,~~
c_m|\l_1,\l_2;\l_3>=0,~~
{\rm for}~i=1,2~~{\rm
and}~~ 
m>0~~.
\end{eqnarray*}
To obtain the highest weight vectors of $\uq1h$, we impose the
conditions:
\begin{eqnarray*}
&&e_i|\l_1,\l_2;\l_3>=0,~~~i=0,1,~{\rm and}~~
h_j|\l_1,\l_2;\l_3>=\l^j|\l_1,\l_2;\l_3>,~~~
j=0,1,2.
\end{eqnarray*}
\noindent Solving these equations, we obtain the following classification:
\begin{itemize}
\item $(\l_1,\l_2;\l_3)=(\b-\a,-\b;\a)$, where $\a$ and $\b$ are 
arbitrary complex numbers. The weight of this vector is $(\l^0,\l^1,\l^2)$=$
(1+\a,-\a,2\b-\a)$. We have
$|(1+\a)\L_0-\a\L_1+(2\b-\a)\b\L_2>$=$|\b-\a,-\b;\a>$.
\end{itemize}
\noindent According to this classification, let us introduce the Fock  
spaces
\bea 
&&\F_{(\a;\b)}=\oplus_{i\in \Z}
F_{\b-\a+i,\;-\b-i;\;\a+i}~~.\label{Fack-space}
\eea
It can be shown that the bosonized action of $\uq1h$ on $\F_{(\a;\b)}$ 
 is closed :$\uq1h \F_{(\a;\b)}=\F_{(\a;\b)}$. 
Hence each Fack space (\ref{Fack-space}) 
constitutes a $\uq1h$-module. However, these modules are  not 
irreducible in general. To obtain the irreducible representations  , we 
introduce a pair of fermionic currents \cite{Kim97,Yan99}
\begin{eqnarray*}
\eta(z)=\sum_{n\in \Z}\eta_nz^{-n-1}=:e^{c(z)}:~,~~
\xi(z)=\sum_{n\in \Z}\xi_nz^{-n}=:e^{-c(z)}:.
\end{eqnarray*}
\noindent The mode expansion of $\eta(z)$ and $\xi(z)$ is well 
defined on $\F_{(\a;\b)}$ with
$\a\in \Z$, and 
the modes satisfy the relations
\bea
&&\xi_m\xi_n+\xi_n\xi_m=\eta_m\eta_n+\eta_n\eta_m=0~,~~
\xi_m\eta_n+\eta_n\xi_m=\delta_{m+n,0}~~.
\eea
Therefore, $\eta_0\xi_0$ and $\xi_0\eta_0$ qualify as
projectors . So we use them to decompose $\F_{(\a;\b)}$ into a 
direct sum $\F_{(\a;\b)}=\eta_0\xi_0\F_{(\a;\b)}\oplus  
\xi_0\eta_0\F_{(\a;\b)}$. Following \cite{Kim97}, 
$\eta_0\xi_0\F_{(\a;\b)}$ is referred to as $Ker_{\eta_0}$ and 
 $\xi_0\eta_0\F_{(\a;\b)}=\F_{(\a;\b)}/\eta_0\xi_0\F_{(\a;\b)}$ as 
$Coker_{\eta_0}$. Since $\eta_0$ commutes (anti-commutes) with 
$\uq1h$, $Ker_{\eta_0}$ and $Coker_{\eta_0}$ are both 
$\uq1h$-modules. 
\subsection{Characters and supercharacters}

In this subsection, we study the character and supercharacter formulas of 
these $\uq1h$-modules which are constructed in the bosonic Fock spaces.
We first of all bosonize the derivation operator $d$ as 
\begin{eqnarray*}
d=-\sum_{1\leq m}\frac{m^2}{[m]_q[2m]_q}\{
A^1_{-m}A^2_m+A^2_{-m}A^1_m+\frac{[2m]_q}{[m]_q}c_{-m}c_m\}
-\frac{1}{2}\{A^1_0A^2_0+c_0(c_0+1)\}.
\end{eqnarray*}
\noindent One can  check that this $d$ obeys the commutation
relations 
\begin{eqnarray*}
[d,h_i]=0,~~~~[d,h^i_m]=mh^i_m,~~~~[d,X^{\pm}_m]=mX^{\pm}_m,~~~
i=1,2,
\end{eqnarray*}
\noindent as required. Moreover, we have $[d,\xi_0]=[d,\eta_0]=0$.

The character and supercharacter of a $\uq1h$-module $M$ are defined by 
\bea
&&Ch_{M}(q,x,y)=tr_M(q^{-d}x^{h_1}y^{h_2})~,\\
&&Sch_{M}(q,x,y)=Str_M(q^{-d}x^{h_1}y^{h_2})=
tr_M((-1)^{N_f}q^{-d}x^{h_1}y^{h_2}),
\eea
\noindent respectively. The Fermi-number operator $N_f$ can be also 
bosonized by $N_f=c_0$. Indeed, such a bosonized operator satisfies 
\begin{eqnarray*}
(-1)^{N_f}\Xi(z)=(-1)^{[\Xi(z)]}\Xi(z),~{\rm for}~ \Xi(z)=X^{\pm,i}(z),
\phi_i(z),\phi^*_i,\psi_i(z),\psi^*_i(z),
\end{eqnarray*}
\noindent as required. Then we have the following result:
\begin{itemize}
\item  (I) {\it Character of $\F_{(\a;\b)}$ for $\a\not\in \Z$}. Since
$\eta_0$ is not defined on this module, it is expected 
that $\F_{(\a;\b)}$ is actually irreducible. We thus have 
\end{itemize}

\begin{Proposition}:
The character and supercharacter of $\F_{(\a;\b)}$ are  
\bea
&&Ch_{\F_{(\a;\b)}}(q,x,y)
=\frac{q^{\frac{1}{2}\a(2\a-2\b+1)}}{\prod_{n=1}^{\infty}(1-q^n)^3}
\sum_{i\in \Z}q^{\frac{1}{2}(i^2+i)}
x^{-\a}y^{2\b-\a+2i},\no\\
&&Sch_{\F_{(\a;\b)}}(q,x,y)
=\frac{q^{\frac{1}{2}\a(2\a-2\b+1)}}{\prod_{n=1}^{\infty}(1-q^n)^3}     
\sum_{i\in \Z}(-1)^{\a+i}q^{\frac{1}{2}(i^2+i)}    
x^{-\a}y^{2\b-\a+2i}.
\eea
\end{Proposition}

\begin{itemize}
\item (II) {\it Characters and supercharacters of $Ker_{\F_{(\a;\b)}}$ 
and  $Coker_{\F_{(\a;\b)}}$,
 for $ \a\in \Z$}. In this case, $\eta_0$ is well defined on
$Ker_{\F_{(\a;\b)}}$ 
and  $Coker_{\F_{(\a;\b)}}$. So we compute the characters and
supercharacters  of these modules 
by using the BRST resolution \cite{Yan99}.
\end{itemize}

Let us define the Fock spaces, for $l\in\Z$
\begin{eqnarray*}
\F^{(l)}_{(\a;\b)}=\oplus_{i\in \Z}F_{\b-\a+i,\;-\b-i;\;\a+i+l}.
\end{eqnarray*}
\noindent We have $\F^{(0)}_{(\a;\b)}=\F_{(\a;\b)}$.
It can be shown that  
$\eta_0$  intertwine these Fock spaces as follows:
\begin{eqnarray*}
&&\eta_0:~\F^{(l)}_{(\a;\b)}\longrightarrow 
\F^{(l+1)}_{(\a;\b)}~,~~\\
&&\xi_0:~\F^{(l)}_{(\a;\b)}\longrightarrow 
\F^{(l-1)}_{(\a;\b)}~.
\end{eqnarray*}
\noindent We have the following  BRST complexes:
\bea
\begin{array}{ccccccc}
\cdots &\stackrel{Q_{l-1}=\eta_0}{\longrightarrow}&
\F^{(l)}_{(\a;\b)}&\stackrel{Q_{l}=\eta_0}{\longrightarrow}&
\F^{(l+1)}_{(\a;\b)}&\stackrel{Q_{l+1}=\eta_0}{\longrightarrow}&
\cdots\\
&&|{\bf O}&&|{\bf O}&&\\
\cdots &\stackrel{Q_{l-1}=\eta_0}{\longrightarrow}&
\F^{(l)}_{(\a;\b)}&\stackrel{Q_{l}=\eta_0}{\longrightarrow}&
\F^{(l+1)}_{(\a;\b)}&\stackrel{Q_{l+1}=\eta_0}{\longrightarrow}&
\cdots
\end{array}\label{BRST}
\eea
\noindent where ${\bf O}$ is an operator such that 
$\F^{(l)}_{(\a;\b)}\longrightarrow \F^{(l)}_{(\a;\b)}$. 
We can get   
\begin{Proposition}\label{brst}
: 
\bea
&&Ker_{Q_l}=Im_{Q_{l-1}},~ ~{\rm for~ any~} l\in \Z,~~~{\rm and }~~
\no\\
&&tr({\bf O})|_{Ker_{Q_l}}=tr({\bf O})|_{Im_{Q_{l-1}}}
=tr({\bf O})|_{Coker_{Q_{l-1}}}~~.
\eea
\end{Proposition}
\noindent {\it Proof}. It follows from the fact that 
$\eta_0\xi_0+\xi_0\eta_0=1$, $(\eta_0)^2=(\xi_0)^2=0$ and 
$\eta_0\xi_0$ ($\xi_0\eta_0$) is the projection operator 
from $\F^{(l)}_{(\a;\b)}$ to $Ker_{Q_{l}}$ 
($Coker_{Q_{l}}$).

\vspace{0.4truecm}

In the following we simply write $Ker_{\eta_0}$
 and $Coker_{\eta_0}$  of $\F_{(\a;\b)}$ as $Ker_{\F_{(\a;\b)}}$ and 
$Coker_{\F_{(\a;\b)}}$, respectively. By proposition \ref{brst}, 
we can compute the characters and 
supercharacters of $Ker_{\F_{(\a;\b)}}$ and  $Coker_{\F_{(\a;\b)}}$, 
for $  \a\in \Z$. We have
%For $\a$
\begin{Proposition}
: The character and
supercharacter of $Ker_{\F_{(\a;\b)}}$ and $Coker_{\F_{(\a;\b)}}$ 
for $\a\in \Z$ are given by
\bea
&&Ch_{Ker_{\F_{(\a;\b)}}}(q,x,y)
=\frac{q^{\frac{1}{2}\a(2\a-2\b+1)}}{\prod_{n=1}^{\infty}(1-q^n)^3}
\sum_{l=1}^{\infty}(-1)^{l+1}q^{\frac{1}{2}(l^2-(1+2\a)l)}
\sum_{i\in
\Z}q^{\frac{1}{2}(i^2+(1-2l)i)}
x^{-\a}y^{2\b-\a+2i},\no\\
&&Ch_{Coker_{\F_{(\a;\b)}}}(q,x,y)
=\frac{q^{\frac{1}{2}\a(2\a-2\b+1)}}{\prod_{n=1}^{\infty}(1-q^n)^3}
\sum_{l=1}^{\infty}(-1)^{l+1}q^{\frac{1}{2}(l^2+(1+2\a)l)}
\sum_{i\in
\Z}q^{\frac{1}{2}(i^2+(1+2l)i)}
x^{-\a}y^{2\b-\a+2i},\no\\
\eea
\noindent and
\bea
&&Sch_{Ker_{\F_{(\a;\b)}}}(q,x,y)
=-\frac{q^{\frac{1}{2}\a(2\a-2\b+1)}}{\prod_{n=1}^{\infty}(1-q^n)^3}
\sum_{l=1}^{\infty}q^{\frac{1}{2}(l^2-(1+2\a)l)}
\sum_{i\in
\Z}(-1)^iq^{\frac{1}{2}(i^2+(1-2l)i)}
(-x)^{-\a}y^{2\b-\a+2i},\no\\
&&Sch_{Coker_{\F_{(\a;\b)}}}(q,x,y)
=-\frac{q^{\frac{1}{2}\a(2\a-2\b+1)}}{\prod_{n=1}^{\infty}(1-q^n)^3}
\sum_{l=1}^{\infty}q^{\frac{1}{2}(l^2+(1+2\a)l)}
\sum_{i\in
\Z}(-1)^iq^{\frac{1}{2}(i^2+(1+2l)i)}
(-x)^{-\a}y^{2\b-\a+2i}.\no\\
\eea
\end{Proposition}
\noindent{\it Proof}. 
Thanks to the resolution of BRST complexes (\ref{BRST}) , 
the trace over $Ker$ and $Coker$ can be 
written as the sum of trace over $\F_{(\a;\b)}^{(l)}$. 
The latter can be 
computed by the technique introduced in \cite{Cla73}.

Note that $\F_{(\a;\b)}^{(1)}=\F_{(\a;\b-1)}$, we have 
\begin{Corollary}
: The following relations hold for any $\a\in \Z$ and $\b$,
\bea
&&Ch_{Coker_{\F_{(\a;\b+1)}}}=Ch_{Ker_{\F_{(\a,\b)}}}~,\\
&&Sch_{Coker_{\F_{(\a;\b+1)}}}=Sch_{Ker_{\F_{(\a,\b)}}}~.
\eea
\end{Corollary}

\subsection{$\uq1h$-module structure of $\F_{(\a;\b+\frac{3}{2}\a)}$}
It is easy to see  that the vector 
\begin{eqnarray*}
|\b+\frac{\a}{2},-\b-\frac{3}{2}\a;\a>\in \F_{(\a;\b+\frac{3}{2}\a)}
\end{eqnarray*}
\noindent plays the role of the highest weight vectors of $\uq1h$-modules.
We  can also check that
\bea 
&&\eta_0|\b+\frac{\a}{2},-\b-\frac{3}{2}\a;\a>=0,~~~{\rm
for}~\a=0,1,2,3,\cdots,\label{Ker1}\\
&&\eta_0|\b+\frac{\a}{2},-\b-\frac{3}{2}\a;\a>\neq 0,~~~{\rm
for}~\a=-1,-2,-3,\cdots.\label{Ker2}
\eea
\noindent It follows that the modules 
$Ker_{\F_{(\a;\b+\frac{3}{2}\a)}}~~(\a=0,1,2,3,\cdots)$  and 
$Coker_{\F_{(\a;\b+\frac{3}{2}\a)}}~~(\a=-1,-2,-3,\cdots)$  are highest 
weight $\uq1h$-modules. Set 
\bea
\l_{\a,\b}=\left\{
\begin{array}{ll}
(1+\a)\L_0-\a\L_1+(2\b+2\a)\L_2& {\rm for }~~\a\not\in\Z\\
(1+\a)\L_0-\a\L_1+(2\b+2\a)\L_2& {\rm for }~~\a=0,1,2,3,\cdots\\
(1+\a)\L_0-\a\L_1+(2\b+2\a+2)\L_2& {\rm for }~~\a=-1,-2,-3,\cdots
\end{array}
\right..
\eea
\noindent Denote by $\V(\l_{\a,\b})$ the highest weight
$\uq1h$-modules with the highest weights $\l_{\a,\b}$. From
(\ref{Ker1})-(\ref{Ker2}) and corollary 1, we obtain 

\begin{Theorem}:
 We have the following identifications of the highest 
weight $\uq1h$-modules:
\bea 
\V(\l_{\a;\b})
&\cong& Ker_{\F_{(\a;\b+\frac{3}{2}\a) }}\equiv 
 Coker_{\F_{(\a;\b+\frac{3}{2}\a+1) }},~~ {\rm for}~~\a\in\Z\no\\
&\cong&\F_{(\a;\b+\frac{3}{2}\a) },~~~~{\rm for}~~~\a\not\in\Z 
\eea
\noindent and when $\a\in\Z$ each Fock space 
$\F_{(\a;\b+\frac{3}{2}\a) }$
can also be 
decomposed explicitly into a direct sum of the highest weight 
$\uq1h$-modules:
\bea 
\F_{(\a;\b+\frac{3}{2}\a) }=\V(\l_{\a,\b})\oplus \V(\l_{\a,\b-1}),~~
  {\rm for}~~\a\in\Z.
\eea  
\end{Theorem}

It is expected that the modules  $\V(\l_{\a,\b})$ are also irreducible with 
respect to the action of $\uq1h$. We thus state

\vspace{0.4truecm}
\noindent{\bf Conjecture 1}: {\it
$\V(\l_{\a,\b})$ are the 
irreducible highest weight $\uq1h$-modules with the highest weight 
 $\l_{\a,\b}$ , i.e 
\bea
&&\V(\l_{\a,\b})=V(\l_{\a,\b}),
\eea
\noindent where $V(\l)$ denotes the irreducible highest weight
$\uq1h$-module with the highest weight $\l$.}

\subsection{Vertex operators  over the irreducible highest
weight $\uq1h$-modules }
In this subsection we study the action of  type I and type II vertex
operators of $\uq1h$  on the irreducible highest weight $\uq1h$-modules.

Using the bosonic representations of the vertex operators
 (\ref{Vertex-operator}), we have the homomorphisms of 
$\uqg2h$-modules:
\bea
&&\phi(z):~\F_{(\a;\b)}\longrightarrow\F_{(\a+1;\b)}\otimes
V_z~,~~
\psi(z):~\F_{(\a;\b)}\longrightarrow  
V_z\otimes\F_{(\a+1;\b)}~,
\label{Inertwining-relation1}\\
&&\phi^*(z):~\F_{(\a;\b)}\longrightarrow\F_{(\a-1;\b
)}\otimes V^{*S}_z~,~~
\psi^*(z):~\F_{(\a;\b)}\longrightarrow
V^{*S}_z\otimes\F_{(\a-1;\b)}~.
\label{Intertwining-relation2}
\eea
We then  consider the vertex operators which intertwine the 
irreducible highest 
weight $\uq1h$-modules. By  conjecture 1, 
we can make  the fellowing identifications:
\bea 
&&\Phi_i(z)=
\left\{
\begin{array}{ll}
\phi_i(z)&{\rm for }~~\a\not\in\Z\\
\eta_0\xi_0\phi_i(z)\eta_0\xi_0& {\rm for}~~\a\in\Z
\end{array}
\right.,~~~
\Phi^*_i(z)=
\left\{
\begin{array}{ll}
\phi^*_i(z)&{\rm for }~~\a\not\in\Z\\
\eta_0\xi_0\phi^*_i(z)\eta_0\xi_0& {\rm for}~~\a\in\Z
\end{array}
\right.,\label{V1}
\no\\
&&\Psi_i(z)=
\left\{
\begin{array}{ll}
\psi_i(z)&{\rm for }~~\a\not\in\Z\\
\eta_0\xi_0\psi_i(z)\eta_0\xi_0& {\rm for}~~\a\in\Z
\end{array}
\right.,~~~
\Psi^*_i(z)=
\left\{
\begin{array}{ll}
\psi^*_i(z)&{\rm for }~~\a\not\in\Z\\
\eta_0\xi_0\psi^*_i(z)\eta_0\xi_0& {\rm for}~~\a\in\Z
\end{array}
\right..\label{V2}\no\\
\eea
\noindent This implies  that  the following vertex operators 
associated with the level-one irreducible highest weight 
$\uq1h$-modules exist:
\bea 
&&\Phi(z)^{\l_{\a+1,\b-3/2}~V}_{\l_{\a,\b}}(z):~~~~
V(\l_{\a,\b})\longrightarrow V(\l_{\a+1,\b-3/2})\otimes
V_z,\no\\
&&\Psi(z)_{\l_{\a,\b}}^{V~\l_{\a+1;\b-3/2}}:~~~~V(\l_{\a,\b})\longrightarrow
V_z\otimes V(\l_{\a+1,\b-3/2}),\no\\
&&\Phi(z)^{\l_{\a-1,\b+3/2}~V^*}_{\l_{\a,\b}}(z):~~~~
V(\l_{\a,\b})\longrightarrow V(\l_{\a-1,\b+3/2})\otimes
V^{*S}_z,\no\\
&&\Psi(z)_{\l_{\a,\b}}^{V^*~\l_{\a-1;\b+3/2}}:~~~~V(\l_{\a,\b})\longrightarrow
V^{*S}_z\otimes V(\l_{\a-1,\b+3/2}).
\eea
\noindent Moreover, the vertex operators defined by 
(\ref{V2}) satisfy the graded Faddeev-Zamolodchikov 
algbera (\ref{ZF1})-(\ref{ZF2}):
\bea
&&\Phi_j(z_2)\Phi_i(z_1)=
\sum_{k,l=1}^{2}R^{kl}_{ij}(\frac{z_1}{z_2})
\Phi_k(z_1)\Phi_l(z_2)(-1)^{[i][j]}~~,\label{ZF3}\\
&&\Psi^*_i(z_1)\Psi^*_j(z_2)=
\sum_{k,l=1}^{2}R^{ij}_{kl}(\frac{z_1}{z_2})
\Psi^*_l(z_2)\Psi^*_k(z_1)(-1)^{[i][j]}~~,\\
&&\Psi^*_i(z_1)\Phi_j(z_2)=
\Phi_j(z_2)\Psi^*_i(z_1)(-1)^{[i][j]},\label{ZF4}
\eea
\noindent and the following invertibility relations:
\bea
\Phi_i(z)\Phi^*_j(z)|_{V(\l_{\a,\b})}=-q\d_{ij}~id|
_{V(\l_{\a,\b})}.
\eea

We can also generalize Miki's construction to the $\uq1h$ case. Define
\begin{eqnarray*}
& &L^+(z)^j_i=\Phi_i(zq^{1/ 2})\Psi^*_j(zq^{-1/2})~~,\\
& &L^-(z)^j_i=\Phi_i(zq^{-1/ 2})\Psi^*_j(zq^{1/2})~~.
\end{eqnarray*}
We have
\begin{Proposition}
: The L-operators $L^{\pm}(z)$ defined above give a realization of the
super Reshetikhin-Semenov-Tian-Shansky algebra \cite{Zha97} at level one  for
the quantum affine superalgebra $\uq1h$ 
\begin{eqnarray*}
& &R(z/w)L^{\pm}_1(z)L^{\pm}_2(w)=L^{\pm}_2(w)L^{\pm}_1(z)R(z/w),\\
& &R(z^{+}/w^{-})L^{+}_1(z)L^{-}_2(w)=L^{-}_2(w)L^{+}_1(z)R(z^-/w^+),
\end{eqnarray*}
\noindent where $L^{\pm}_1(z)=L^{\pm}(z)\otimes 1$,
$L^{\pm}_2(z)=1\otimes L^{\pm}(z)$ and $z^{\pm}=zq^{\pm\frac{1}{2}}$.
\end{Proposition}

\noindent {\it Proof}. Straightforward by using (\ref{ZF3})-(\ref{ZF4}).

\vskip.3in
\subsection*{ Acknowledgements.}

This work has been financially supported by the Australian Research 
Council large, small and QEII fellowship grants. We would like to
thank Ryu Sasaki for a careful reading of the manuscript and many
helpful comments. We thank Bo-Yu Hou for
encouragement and useful discussions. W.L.Yang thanks
Y.-Z.Zhang and department of Mathematics, University of Queensland,
for their kind  hospitality. W.L.Yang was also partially supported by the
National Natural Science Foundation of China.

%\newpage

\end{document}